\font\piccolo=cmr8.
\font\script=eusm10.
\font\sets=msbm10.
\font\stampatello=cmcsc10.
\font\symbols=msam10.

\def\1{{\bf 1}}
\def\uH{u_{\!_H}\!}
\def\UH{U_{\!_H}\!}
\def\modUH{\widetilde{U}_{\!_H}\!}
\def\hatuH{\widehat{u}_{\!_H}\!}	
\def\hatUH{\widehat{U}_{\!_H}\!}

\def\avesum{\sum_{x\sim N}}
\def\square{\hbox{\vrule\vbox{\hrule\phantom{s}\hrule}\vrule}}
\def\defineq{\buildrel{def}\over{=}}

\def\doublesum{\mathop{\sum\sum}}
\def\integrale{\mathop{\int}}

\def\C{\hbox{\sets C}}
\def\D{\hbox{\sets D}}
\def\N{\hbox{\sets N}}

\def\R{\hbox{\sets R}}

\def\Corr{\hbox{\script C}\!}
\def\EssBdd{\hbox{\symbols n}\,}
\def\modSel{{\widetilde{J}}}

\par
\centerline{\bf Some optimal links between generations of correlation averages}
\bigskip
\centerline{\stampatello giovanni coppola\footnote{$^1$}{\piccolo titolare di un Assegno \lq \lq Ing.Giorgio Schirillo\rq \rq \thinspace dell'Istituto Nazionale di Alta Matematica (Fellow \lq \lq Ing.Giorgio Schirillo\rq \rq \thinspace of the Istituto Nazionale di Alta Matematica).} - maurizio laporta}
\bigskip
{

{\par
\noindent
{\bf Abstract.} For a real-valued and essentially bounded arithmetic function $f$, i.e., $f(n)\ll_{\varepsilon}\!n^{\varepsilon},\,\forall\varepsilon\!>\!0$, we \enspace give some optimal links between non-trivial  bounds for the sums 
$\sum_{h\le H}\sum_{N<n\le 2N}f(n)f(n-h)$, $\sum_{N<x\le 2N} \big| \sum_{x<n\le x+H}f(n)\big|^2$ and
$\sum_{N<n\le 2N}  \big| \sum_{0\le |n-x|\le H}\big(1-{|n-x|\over H}\big)f(n)\big|^2$, with $H=o(N)$ as $N\to\infty$.

\footnote{}{\par \noindent {\it Mathematics Subject Classification} $(2010) : 11{\rm N}37
.$}
}
\bigskip

\par
\noindent {\bf 1. Introduction.} 
\bigskip
\par
\noindent
The {\it correlation} of a complex-valued arithmetic function $f$
is a shifted convolution sum of the form
$$
\Corr_f(h)\defineq \sum_{n\sim N}f(n)\overline{f(n-h)}, 
$$
\par
\noindent
where $n\sim N$ means that $n\in(N,2N]\cap\N$, while $N$ and the {\it shift} $h$ are integers such that $|h|\le N$.
In particular, this allows restricting $f$ to $1\le n\le 3N$ when dealing with $\Corr_f(h)$. 
Further,  note that 
$$
\Corr_f(h)
=\doublesum_{{n\sim N \thinspace m\sim N}\atop {n-m=h}}f(n)\overline{f(m)}+O(\Vert f\Vert_{\infty}^2|h|),\quad
\hbox{with}\ \Vert f\Vert_{\infty}\defineq\max_{1\le n\le 3N}|f(n)|.
$$
\par
\noindent
In [CL1] and [CL3] we have investigated the 
connection between the correlations of $f$ and its {\it Selberg integral}
$$
J_f(N,H)\defineq \avesum \Big| \sum_{x<n\le x+H}f(n)-M_f(x,H)\Big|^2\, , 
$$
\par
\noindent
where $M_f(x,H)$ is the expected {\it mean value} of $f$ in short intervals, with $H=o(N)$ as $N\to\infty$
(to avoid trivialities, hereafter we assume that $H\to\infty$). 
More in general, we extended such an investigation to weighted versions of $J_f(N,H)$, that include, also, 
the so-called {\it modified Selberg integral}
$$
\modSel_f(N,H)\defineq \avesum \Big| \sum_{n}C_H(n-x)f(n)-M_f(x,H)\Big|^2,
$$
\par
\noindent
where the {\it Cesaro weight} 
$C_H(t)\defineq\max(1-|t|/H,0)$ allows taking the same mean value that appears in $J_f(N,H)$.
Indeed, in [CL1] (see Lemma 7 there), by means of an elementary Dispersion Method, 
it is shown that weighted Selberg integrals, for a wide class of arithmetic functions, are actually linked to averages of their correlations (see also the proof of Lemma 1 below). Mainly inspired by the prototype of the divisor function $d_k$, 
we sticked to the case of a real-valued and {\it essentially bounded} $f$, i.e. $f(n)\ll_{\varepsilon} n^{\varepsilon}$ ($\forall\varepsilon>0$). Here we recall that $\ll$ is Vinogradov's notation, synonimous to Landau's $O$-notation. In particular, $\ll_{\varepsilon}$ means that the implicit constant might depend on an arbitrarily small $\varepsilon>0$, which might change at each occurrence. Moreover,
we abbreviate
$A(N,H)\EssBdd B(N,H)$
whenever $A(N,H)\ll_{\varepsilon} N^{\varepsilon}B(N,H)$, $\forall \varepsilon>0$.\par

In [CL1] we also searched for links between non-trivial bounds for  $J_f(N,H)$, $\modSel_f(N,H)$ and the so-called {\it deviation} of $f$, that is defined as
$$
\D_f(N,H)\defineq\sum_{h\le H}\Corr_f(h)-\avesum {M_f(x,H)^2\over H}.
$$
In this paper we focus on the special case of a real-valued and essentially bounded function $f$ which is also {\it balanced}, 
that is $M_f(x,H)$ vanishes identically. Therefore, this yields 
$$
\eqalign
{
\D_f(N,H)=&\sum_{h\le H}\Corr_f(h)\, ,\cr
J_f(N,H)=&\avesum \Big| \sum_{x<n\le x+H}f(n)\Big|^2\, ,\cr
\modSel_f(N,H)=&\avesum \Big| \sum_{n}C_H(n-x)f(n)\Big|^2\, .\cr
}
$$

\noindent 
In this regard, we have the following result that improves on the bounds given in [CL1]. 
\bigskip

\noindent 
{\stampatello Theorem.} 
{\it Let $H,N\in\N$ with $H=o(N)$ as $N\to\infty$, and let $f:\N \rightarrow \R$ be essentially bounded and balanced. 
\par
\noindent
i) If there is $A\in [-1,1)$ such that $J_f(N,H)\EssBdd NH^{1+A}$ and 
$J_f(N,H_1)\EssBdd NH_1^{1+A}$ 
for $H_1=\big[H^{1-{{2(1-A)}\over {3-A}}}\big]$, then
$$
\D_f(N,H)\EssBdd (N+H^{2-A})H^{1-{{2(1-A)}\over {3-A}}}\, . 
$$
\par
\noindent
ii) If there is $A\in [-3,1)$ such that $\modSel_f(N,H)\EssBdd NH^{1+A}$ and 
$\modSel_f(N,H_2)\EssBdd NH_2^{1+A}$ 
for $H_2=\big[H^{1-{{2(1-A)}\over {5-A}}}\big]$, then
$$
J_f(N,H)\EssBdd (N+H^{2-A})H^{2-{{4(1-A)}\over {5-A}}}\, . 
$$
}
\par
\noindent
Hereafter, $[t]$ denotes the {\it integer part} of $t\in\R$.
\bigskip

\noindent
{\stampatello Remarks.}\par\noindent  
1. Note that 
$$
A\in [-1,1)\Leftrightarrow 1+A\in [0,2)\Leftrightarrow 1-{{2(1-A)}\over {3-A}}\in[0,1)\Leftrightarrow (2-A)^{-1}\in [1/3,1),
$$
$$
A\in [-3,1)\Leftrightarrow 1+A\in [-2,2)\Leftrightarrow 1-{{2(1-A)}\over {5-A}}\in[0,1)\Leftrightarrow (2-A)^{-1}\in [1/5,1).
$$
\par
\noindent 
Finally, the minimal values $-1$ and $-3$ for $A$, resp., in $i)$ and $ii)$, are sharp. To see this, it is enough to consider
the function $f(n)=(-1)^{n+1}$, which is essentially bounded and  balanced  (its Dirichlet series is known as {\it Dirichlet's eta function}). Since its correlation is $\Corr_f(h)=(-1)^{h}N$, then 
$$
\eqalign
{
\D_f(N,H)=&N\sum_{h\le H}(-1)^{h}=\cases{-N\ &\hbox{if}\ $H$\ \hbox{is odd}\cr 0\ &\hbox{otherwise},\cr}
\cr
J_f(N,H)=&\avesum \Big| \sum_{x<n\le x+H}(-1)^n\Big|^2=\cases{N\ &\hbox{if}\ $H$\ \hbox{is odd}\cr 0\ &\hbox{otherwise},\cr}
\cr
\modSel_f(N,H)=&{1\over H^2}\avesum \Big|\sum_{h\le H}\sum_{0\le |n-x|<h}(-1)^n\Big|^2\cr
                =&{1\over H^2}\avesum \Big|\sum_{h\le H\atop x+h\, {\rm odd}}1-\sum_{h\le H\atop x+h\, {\rm even}}1\Big|^2
                 =\cases{N/H^2\ &\hbox{if}\ $H$\ \hbox{is odd}\cr 0\ &\hbox{otherwise},\cr}\cr
}
$$
where we have used the Cesaro identity 
$$
\sum_{0\le |n-x|\le H}\left( 1-{{|n-x|}\over H}\right)f(n)={1\over H}\sum_{h\le H}\sum_{0\le |n-x|<h}f(n)\ .
$$

\noindent
2. Through an application of the so-called {\it length inertia} (see [CL3] and [CL4]),  it could be shown that the Theorem hypothesis on  $J_f(N,H)$ and $\modSel_f(N,H)$ might be dropped  without affecting the result. 
\vfill
\eject

\par				
\noindent {\bf 2. Lemmata.} 
\bigskip
\par
\noindent
First, let us introduce some notation and some auxiliary functions. For the {\it unit step} weight
$$
\uH(a)\defineq \cases{1\ &\hbox{if}\ $a\in[1,H]\cap\N$\cr 0\ &\hbox{otherwise},\cr}
$$
\par
\noindent
we set
$$
\eqalign{
&\UH(h)\defineq \doublesum_{{a \thinspace \quad \thinspace b}\atop {b-a=h}}\uH(b)\uH(a)\, ,\quad
\modUH(h)\defineq{1\over H^2}\doublesum_{{a \thinspace \quad \thinspace b}\atop{b-a=h}}\UH(b)\UH(a)\, ,
\cr
&\hatuH(\beta)\defineq \sum_{h\le H}e(h\beta)\, ,\qquad\qquad\ \hatUH(\beta)\defineq \sum_{h\le H}\UH(h)e(h\beta),\quad
\forall \beta\in \R.\cr
}
$$
where $e(\alpha)\defineq e^{2\pi i\alpha}$.
Moreover, for the function $f$ under consideration we denote
$$
\widehat{f}(\beta)\defineq \sum_{n\sim N}f(n)e(n\beta)\, .
$$
Next Lemma is a consequence of Lemma 7 from [CL1]. Somehow the formul\ae\ $(I)-(III)$ were already 
implicit between the lines of [CL1], where, however, the underlying assumption  that $f$ has to be also essentially bounded is, in fact, redundant.
\bigskip

\par
\noindent
{\stampatello Lemma 1.} {\it For every balanced $f:\N\rightarrow\R$ one has
$$
\leqalignno
{
\D_f(N,H)=& \int_{0}^{1}\big| \widehat{f}(\beta)\big|^2\, \hatuH(-\beta){\rm d}\beta+ O(H^2 \Vert f\Vert_{\infty}^2),&
(I)\cr
J_f(N,H)
=&\int_{0}^{1}\big| \widehat{f}(\beta)\big|^2 \big|\hatuH(\beta)\big|^2 {\rm d}\beta +O(H^3 \Vert f\Vert_{\infty}^2),
& (II)\cr
\modSel_f(N,H)
=&
\int_{0}^{1}\big| \widehat{f}(\beta)\big|^2 {\big|\hatuH(\beta)\big|^4\over H^2}
{\rm d}\beta + O(H^3 \Vert f\Vert_{\infty}^2).
&(III)\cr
}
$$
}
\bigskip

\par
\noindent
{\stampatello Proof.} Since 
$$
\D_f(N,H)=\sum_{h\le H}\Corr_f(h),
$$
then $(I)$ follows immediately because it is plain that
$$
\Corr_f(h)=\doublesum_{{n\sim N \thinspace m\sim N}\atop {n-m=h}}f(n)f(m)+O(\Vert f\Vert_{\infty}^2|h|)
=\int_{0}^{1}\big| \widehat{f}(\beta)\big|^2 e(-h\beta){\rm d}\beta + O(\Vert f\Vert_{\infty}^2 |h|).
$$
In order to show $(II)$ and $(III)$ let us recall that the Selberg integral and the modified one for any real arithmetic function $f$ are related to the correlation averages (see [CL1], Lemma 7), respectively as
$$ 
\eqalign{
J_f(N,H)=&\sum_h \UH(h)\Corr_f(h)-2\sum_n f(n)\avesum \uH(n-x)M_f(x,H)+\avesum M_f(x,H)^2
+ O(H^3 \Vert f\Vert_{\infty}^2), \cr
\modSel_f(N,H)=&\sum_h \modUH(h)\Corr_f(h)-{2\over H}\sum_n f(n)\avesum \UH(n-x)M_f(x,H)+\avesum M_f(x,H)^2+ O(H^3 \Vert f\Vert_{\infty}^2). \cr
}
$$
\par
\noindent
In particular, by setting $M_f(x,H)=0$ in these formul\ae\ and by using the properties
$$
\UH(h)
=\sum_{a\le H-|h|}1\ ,\qquad 
\hatUH(\beta)=|\hatuH(\beta)|^2,
$$
\par				
\noindent
we get $(II)$ and $(III)$, because it is easily seen that
$$
\sum_{h}\UH(h)\Corr_f(h)=\sum_{h\le H}\sum_{0\le |a|<h}\Corr_f(a)=\int_{0}^{1}\big| \widehat{f}(\beta)\big|^2 \hatUH(-\beta){\rm d}\beta + O(H^3 \Vert f\Vert_{\infty}^2),
$$
$$
\sum_{h}\modUH(h)\Corr_f(h)=\sum_{h}\doublesum_{{a \thinspace \quad \thinspace b}\atop{b-a=h}}{{\UH(b)\UH(a)}\over H^2}\Corr_f(h)={1\over H^2}\int_{0}^{1}\big| \widehat{f}(\beta)\big|^2 {\big|\hatUH(\beta)\big|^2}{\rm d}\beta + O(H^3 \Vert f\Vert_{\infty}^2).
$$
The Lemma is completely proved.\hfill $\square$ 

\bigskip

\par
\noindent {\stampatello Remark.} Consistently with the terminology introduced in [CL1], we refer to $(I)$, $(II)$ and $(III)$ as 
a {\it first, second} and {\it third generation} formula, respectively. As transpires also from the above proof, such formul\ae\ correspond to iterations of correlations' averages.
\smallskip

Next Lemma gives two versions of a Gallagher's inequality (see [Ga], Lemma 1), that have been discussed in [CL2] and [CL4]. 
\bigskip
\par
\noindent {\stampatello Lemma 2.}
{\it Let $N,h\in\N$ be such that $h\to \infty$ and $h=o(N)$ as $N\to \infty$. If \thinspace $f:\N \rightarrow \C$ \thinspace is essentially bounded and balanced,
then}
$$
h^2\int_{-{1\over {2h}}}^{{1\over {2h}}}| \widehat{f}(\alpha)|^2 {\rm d}\alpha\,
\EssBdd J_f(N,h)+h^3, 
\leqno{1)}
$$
$$
h^2\int_{-{1\over {2h}}}^{{1\over {2h}}}| \widehat{f}(\alpha)|^2 {\rm d}\alpha\, 
\EssBdd \modSel_f(N,h)+h^3. 
\leqno{2)}
$$
}

\par
\noindent {\bf 3. Proof of the Theorem.} 
\bigskip
\par
\noindent
In what follows, we will appeal to the well-known property
$$
|\hatuH(\alpha)|
={{|\sin(\pi H \alpha)|}\over {\sin(\pi \alpha)}}\le
{{1}\over { \sin(\pi \alpha)}}< {{1}\over {2\alpha}},\quad\hbox{for}\ 0<\alpha<{1\over 2}.
$$
In particular, this yields the implication
$$
|\hatuH(\alpha)|
>h
\enspace \Longrightarrow \enspace
|\alpha|< {1\over {2h}}\, .
\leqno{(\ast)}
$$
\bigskip
\par
\noindent {\stampatello Proof of $i)$.} Since $f$ is essentially bounded and balanced, then from $(I)$ of Lemma 1 we infer
$$
\eqalign{
\D_f(N,H)\EssBdd& \int_{-1/2}^{1/2}\left|\widehat{f}(\alpha)\right|^2 \left|\hatuH(\alpha)\right|{\rm d}\alpha + H^2\cr
\EssBdd&
H^{1-\delta} \integrale_{\left|\hatuH(\alpha)\right|\le [H^{1-\delta}]}\left|\widehat{f}(\alpha)\right|^2 {\rm d}\alpha
 + H^{1-\gamma} \integrale_{[H^{1-\delta}]<\left|\hatuH(\alpha)\right|\le H^{1-\gamma}}\left|\widehat{f}(\alpha)\right|^2 {\rm d}\alpha\cr
 &
 + H^{\gamma-1} \integrale_{\left|\hatuH(\alpha)\right|>H^{1-\gamma}}\left|\widehat{f}(\alpha)\right|^2 |\hatuH(\alpha)|^2 {\rm d}\alpha + H^2, \cr
 }
$$
\par				
\noindent
where $\gamma,\delta$ are real numbers to be determined later, so that
$0<\gamma \le \delta$. Thus, 
by applying $(\ast)$, Parseval's identity and $(II)$ of Lemma 1, we get
$$
\D_f(N,H)\EssBdd NH^{1-\delta} + H^{1-\gamma} \integrale_{|\alpha|\le {1\over {2H_1}}}\left|\widehat{f}(\alpha)\right|^2 {\rm d}\alpha
 + J_f(N,H)H^{\gamma-1} + H^{2+\gamma}, 
$$
\par
\noindent
where we have set $H_1=[H^{1-\delta}]$. Thus, by $1)$ of Lemma 2 and, then, assuming that \enspace $J_f(N,H_1)\EssBdd NH_1^{1+A}$ \enspace 
and \enspace $J_f(N,H)\EssBdd NH^{1+A}$, 
we write
$$
\eqalign{
\D_f(N,H)\EssBdd& NH^{1-\delta} + H^{1-\gamma}H_1^{-2} J_f(N,H_1) + H^{1-\gamma}H_1+ J_f(N,H)H^{\gamma-1}
+ H^{2+\gamma} \cr
\EssBdd&NH(H^{-\delta}+ H^{-\gamma-(1-\delta)(1-A)}+H^{\gamma+A-1})+H^{2+\gamma}. \cr
}
$$
\par
\noindent
Now, observe that $\delta=\gamma+(1-\delta)(1-A)=1-A-\gamma$ is satisfied by 
$\delta={{2(1-A)}\over {3-A}}$ and $\gamma={{(1-A)^2}\over {3-A}}$, which obey the condition $0<\gamma \le \delta$ whenever 
$A\in[-1,1)$. This yields the inequality for $\D_f(N,H)$ stated in $i)$.
\smallskip

\par
\noindent {\stampatello Proof of $ii)$.} Since we closely follow the proof of $i)$, we skip some details. By using $(II)$ and $(III)$ of Lemma 1 and applying $2)$ of Lemma 2, as before we can write
$$
\eqalign{
J_f(N,H)\EssBdd& NH^{2-2\delta} + H^{2-2\gamma} \integrale_{|\alpha|\le {1\over {2H_2}}}\left|\widehat{f}(\alpha)\right|^2 {\rm d}\alpha
 + H^{2\gamma} \int_{-1/2}^{1/2}\left|\widehat{f}(\alpha)\right|^2 {{|\hatuH(\alpha)|^4}\over {H^2}} {\rm d}\alpha + H^3\cr
 \EssBdd& NH^{2-2\delta} + H^{2-2\gamma}H_2^{-2}\modSel_f(N,H_2)+H^{2-2\gamma}H_2 + H^{2\gamma}\modSel_f(N,H)+H^{3+2\gamma}, \cr
 }
$$
\par
\noindent
where $H_2=[H^{1-\delta}]$. Thus, from $\modSel_f(N,H)\EssBdd NH^{1+A}$ and $\modSel_f(N,H_2)\EssBdd NH_2^{1+A}$, it follows
$$
J_f(N,H)\EssBdd NH^2(H^{-2\delta} + H^{-2\gamma-(1-A)(1-\delta)}
 + H^{A-1+2\gamma})+ H^{3+2\gamma}. 
$$
\par
\noindent
The conclusion follows by taking  $\delta={{2(1-A)}\over {5-A}}$ and $\gamma ={{(1-A)^2}\over {2(5-A)}}$, noticing that
$0<\gamma \le \delta$ whenever $A\in[-3,1)$. The Theorem is completely proved.\hfill $\square$

\bigskip

\par
\centerline{\bf References}
\smallskip
\par
\item{\bf [CL1]} G. Coppola, M. Laporta, {\it Generations of correlation averages}, Journal of Numbers, Vol. 2014 (2014), Article ID 140840, 13 pages, http://dx.doi.org/10.1155/2014/140840 (draft at arxiv:1205.1706v3)
\smallskip
\item{\bf [CL2]} G. Coppola, M. Laporta, {\it A modified Gallagher's Lemma},  preprint at  arxiv.org/abs/1301.0008v1
\smallskip
\item{\bf [CL3]} G. Coppola, M. Laporta, {\it Symmetry and short interval mean-squares}, (submitted), preprint available at 
arXiv:1312.5701v1.
\smallskip
\item{\bf [CL4]} G. Coppola, M. Laporta, {\it A generalization of Gallagher's Lemma for exponential sums}, to appear on \v{S}iauliai Mathematical Seminar, (draft at arxiv.org/abs/1411.1739v1)
\smallskip
\item{\bf [Ga]} P. X. Gallagher, {\it A large sieve density estimate near $\sigma =1$}, Invent. Math. {\bf 11} (1970), 329--339. $\underline{\tt MR\enspace 43\# 4775}$ 

\bigskip

\leftline{\tt Giovanni Coppola \hfill Maurizio Laporta}
\leftline{\tt Universit\`a degli Studi di Napoli \hfill Universit\`a degli Studi di Napoli}
\leftline{\tt Home address \negthinspace : \negthinspace Via Partenio \negthinspace 12 \negthinspace - \hfill Dipartimento di Matematica e Appl.}
\leftline{\tt - 83100, Avellino(AV), ITALY \hfill Compl.Monte S.Angelo}
\leftline{\tt e-page : $\! \! \! \! \! \!$ www.giovannicoppola.name \hfill Via Cinthia - 80126, Napoli, ITALY}
\leftline{\tt e-mail : $\! \! \! \! \! \!$ giovanni.coppola@unina.it \hfill e-mail : mlaporta@unina.it}

\bye